\documentclass[12pt]{article}

\usepackage[english]{babel}
\usepackage{amssymb,amsmath}
\usepackage{hyperref}

\newcommand{\footremember}[2]{%
    \footnote{#2}
    \newcounter{#1}
    \setcounter{#1}{\value{footnote}}%
}

\newtheorem{theorem}{Theorem }[section]

\newtheorem{corollary}[theorem]{Corollary}

\newcommand{\qed}{\hspace*{\fill}$\Box$}

\title{ Paley type partial difference sets in abelian groups}

\author{Zeying Wang \footremember{MTU}{Michigan Technological University}  \footnote{zeying@mtu.edu}}
\date{}
\begin{document}

\maketitle 

\abstract{Partial difference sets with parameters $(v,k,\lambda,\mu)=(v, (v-1)/2, (v-5)/4, (v-1)/4)$ are called Paley type partial difference sets. In this note we prove that if  there exists a Paley type partial difference set in an abelian group $G$ of  an order not a prime power, then $|G|=n^4$ or $9n^4$, where  $n>1$ is an odd integer.  In 2010, Polhill 
\cite{Polhill} constructed Paley type partial difference sets in abelian groups with those orders. Thus, combining with the constructions of Polhill and the classical Paley construction using non-zero squares of a finite field,  we completely answer the following question: `` For which odd positive integer $v > 1$, can we find a Paley type partial difference set in an abelian group of order $v$?"}

\section{Introduction and main result}

Let $G$ be a finite abelian group of order $v$, and let $D\subseteq G$ be a subset of size $k$. We say that $D$ is a $(v,k,\lambda,\mu)$-{\it partial difference set} (PDS) in $G$ if the expressions $gh^{-1}$, $g$, $h\in D$, $g\neq h$, represent each non-identity element in $D$ exactly $\lambda$ times, and each non-identity element of $G$ not in $D$ exactly $\mu$ times. If we further assume that $D^{(-1)}=D$  (where $D^{(s)}=\{g^s:g\in D\}$) and  $e \notin D$ (where $e$ is the identity element of $G$), then $D$ is called a {\it regular} partial difference set.  A regular PDS is called {\it trivial} if $D\cup\{e\}$ or $G\setminus D$ is a subgroup of $G$. The condition that $D$ be regular is not a very restrictive one, as  $D^{(-1)}=D$ is automatically fulfilled whenever $\lambda\neq\mu$, and $D$ is a PDS if and only if $D\cup\{e\}$ is a PDS. \medskip

Throughout this paper, we will use the following standard notation: $\beta=\lambda-\mu$ and  $\Delta=\beta^2+4(k-\mu)$.

\medskip

If $v \equiv 1 \pmod 4$, a regular $(v, (v-1)/2, (v-5)/4, (v-1)/4)$-PDS  is called  a {\it Paley type partial difference set}.  In 1933, Paley \cite{Paley} discovered the following well-known family of PDSs.

{\bf Example:} Let  $q$ be a prime power with $q \equiv 1 \pmod 4$.  The set of nonzero squares in the finite field $\mathbb{F}_q$ is a Paley type PDS in the group $(\mathbb{F}_q, +)$.

Athough many further examples  of Paley-type PDSs in $p$-groups were found (for example \cite{Davis}, \cite{Leung}),  not much was known about Paley type PDSs in abelian groups of order not a prime power.  In 1994,  Arasu, Jungnickel, Ma and Pott  \cite{Arasu} asked the following question:

{\bf Question 1:} Suppose $G$ is an abelian group of order $v\equiv 1 \pmod 4$. If $v$ is not a prime power, does there exist a Paley type PDS in $G$?

\medskip

In 2009, Polhill  \cite{Polhill2009} gave the first construction of Paley type PDSs in groups having a non-prime power order (namely in $\mathbb{Z}_3^2\times \mathbb{Z}_p^{4t}$, $p$ an odd prime).  In 2010, Polhill  \cite{Polhill} gave a much more general construction and showed that Paley type PDSs exist in some abelian groups with orders $n^4$ and $9n^4$ for all  $n$  where $n$ is odd and $n>1$.  Thus it is natural to ask the following question:

{\bf Question 2:}  For which odd positive integer $v > 1$, can we find a Paley type PDS in an abelian group of order $v$?

\medskip

Before we get into more details of answering the above question, we will mention  an important result that  was proved by  Ma \cite{MA84} in 1984.

\begin{theorem}{\rm \cite{MA84}} \label{Ma84}
  Let $D$ be an abelian regular $(v, k, \lambda, \mu)$-PDS, and assume that $\Delta=(\lambda-\mu)^2+4(k-\mu)$ is not a perfect square. Then $D$ is of Paley type; more precisely, $D$ has parameters
$$\left(p^{2s+1}, \frac{p^{2s+1}-1}{2}, \frac{p^{2s+1}-5}{4}, \frac{p^{2s+1}-1}{4}\right),$$
where $p$ is a prime congruent to 1 modulo 4.
\end{theorem}

Let  $D$ be a regular Paley type PDS in an abelian group $G$, where $|G|=v$.  Then $\Delta=(-1)^2+4(\frac{v-1}{4})=v$.  If $v$ is not a square,  then by Theorem \ref{Ma84}, we have $|G|=v=p^{2s+1}$ for some prime $p\equiv 1 \pmod 4$. For a  prime power $q\equiv 1 \pmod 4$, we can always construct a Paley type PDS in $(\mathbb{F}_q, +)$ using the non-zero squares of the finite field $\mathbb{F}_q$. Thus to answer Question 2, we only need to focus on the existence of Paley type PDSs   when $v$ is a perfect square and not a prime power, that is, when $v=p_1^{2t_1}p_2^{2t_2}\cdots p_k^{2t_k}$, $k \ge 2$, $p_1$, $p_2$, $\cdots$, $p_k$ are distinct odd prime numbers. In \cite{Wang2019}, we partially answered this question.  In this note, we will completely answer Question 2, and we record our main result below:

\begin{theorem}\label{restriction_order}
Let $v$ be an odd positive integer $>1$. Then there exists a Paley type PDS in some abelian group $G$ of order $v$ if and only if $v$ is a prime power and $v \equiv 1 \pmod 4$, or $v=n^4$ or $9n^4$, with $n>1$ an odd positive integer.
\end{theorem}

\section{Proof of Theorem \ref{restriction_order}}

Below we cite a result on ``sub-partial difference sets"  which was discovered by Ma in \cite{MA94a}, and also recorded in the 1994 survey of Ma \cite{MA94b}.  We quote this result in the form given in \cite{Arasu}. 

\begin{theorem}{\rm \cite{Arasu}}\label{Ma1}
Let $D$ be a non-trivial regular $(v,k,\lambda,\mu)$-PDS in an abelian group $G$, assume that $\Delta$ is a square, say $\Delta=\delta^2$. If $H$ is a subgroup of $G$ such that $(|H|, |G|/|H|)=1$ and $|G|/|H|$ is odd, then $D_1=D\cap H$ is a regular $(v_1,k_1, \lambda_1, \mu_1)$-PDS in $H$ with parameters
$$v_1=|H|,  \quad \beta_1=\beta-2\theta \pi, \quad\Delta_1=\pi^2, $$
$$k_1=\frac{1}{2}\left[(v_1+\beta_1)\pm \sqrt{(v_1+\beta_1)^2-(\Delta_1-\beta_1^2)(v_1-1)}\right],$$

where $\pi=\gcd(|H|, \delta)$ and $\theta$ is the integer satisfying $(2\theta-1)\pi\le \beta <(2\theta+1)\pi$. Moreover, if $D_1\neq \emptyset$, $D_1 \neq H \setminus\{e\}$ and if $\delta=p^r \pi$, where $p\ge 5$ is a prime and $\pi>1$ is relative prime to $p$, then either $r$ is even and $\theta \equiv 0 \pmod{p-1}$ or $r$ is odd and $\theta \equiv (p-1)/2\pmod{p-1}$.
\end{theorem}

{\bf Note:} In Theorem \ref{Ma1}, if $\Delta_1=|H|$, we have $$k_1=\frac{1}{2}\left[ |H| +\beta_1\pm (\beta_1+1)\sqrt{|H|}  \right].$$

\medskip

Now we state our first result:

\begin{theorem}\label{r_condition}
Let $G$ be an abelian group of order $p^{2r}u^2$, $p \ge 5$ a prime number, $u>1$, and $\gcd(p,u)=1$. If $D$ is a Paley type partial difference set in $G$, then $r$ is even.
\end{theorem}

{\bf Proof:} We prove this by contradiction. Assume $r$ is odd. Since $D$ is a Paley type partial difference set in $G$, we have $\Delta=\beta^2+4(k-\mu)=p^{2r}u^2$, $\delta=\sqrt{\Delta}=p^r u$. Let $H$ be a subgroup of $G$ of order $u^2$. By Theorem \ref{Ma1} we have $\pi=\gcd(|H|, \delta)=u$, and $D_1=D\cap H$ is also a regular partial difference set in $H$, and $\theta$ is the integer satisfying $(2\theta-1)\pi\le \beta=-1 <(2\theta+1)\pi$. Thus $\theta=0$, and $\beta_1=\beta-2\theta \pi=-1$. Since $\Delta_1=\pi^2=u^2=|H|$, by Theorem \ref{Ma1} we know that  $$|D_1|=k_1=\frac{u^2-1}{2}.$$ Since $u>1$, we have $D_1 \neq \emptyset$. Also $D\neq H\setminus \{e\}$. As $\pi=u>1$ and $\gcd(p,u)=1$, by Theorem \ref{Ma1}, if $r$ is odd, we have 
$$\theta \equiv (p-1)/2 \pmod{p-1},$$
contradicting with the fact that $\theta=0$.  Thus the conclusion follows.
\qed.

\bigskip

By applying Thorem \ref{r_condition}, we can prove the following Corollary:

\begin{corollary}\label{cor_1}
Let $G$ be an abelian group, and $|G|=p_1^{2t_1}p_2^{2t_2}\cdots p_k^{2t_k}$, $k \ge 2$, $p_1$, $p_2$, $\cdots$, $p_k$ are distinct odd prime numbers.  If $D$ is a Paley type PDS in $G$, then $|G|=n^4$ or $9n^4$, with $n>1$ an odd positive integer. 
\end{corollary} 

{\bf Proof:}  For any $1\le i \le k$, if $p_i \ge 5$, by letting $u=|G|/p_i^{2t_i}$ and applying Theorem \ref{r_condition}, we have that $t_i$ must be an even number. Thus it follows that $|G|=n^4$ or $9n^4$, with $n>1$ an odd positive integer.
\qed.

\medskip

On the other hand, in 2010 J. Polhill  \cite{Polhill} proved the following results:

\begin{theorem}\label{Polhill}
Let $n$ be a positive odd number with $n>1$. Then there is a Paley type partial difference set in a group of order $n^4$ and $9n^4$, respectively.
\end{theorem}

\medskip

Combining Corollary \ref{cor_1} and Theorem \ref{Polhill}, Theorem \ref{restriction_order} follows.

\bigskip

{\bf Acknowledgement:} After the note was written, the author noticed that in \cite{Leung}, Leung and Ma listed a very similar result to Theorem \ref{r_condition} (Lemma 2.4 in \cite{Leung}), but no proof of it was given before.

\end{document}